\documentclass [11pt]{amsart} 
\usepackage{epic,eepic,latexsym, amssymb, amscd, amsfonts, color}

\input xy
\xyoption {all}
\def \gbar {{\bar g}}
\def \s {{\mathcal {S}}}
\def \u {{\mathcal U}}
\def \m {{\mathcal M}}
\def \p {{\mathcal P}_{r,k}}

\newtheorem{theorem}{Theorem} 
\newtheorem {lemma}{Lemma} 
\newtheorem {corollary}{Corollary} 
\newtheorem {proposition}{Proposition}

\newcommand{\comment}[1]{}
\newtheorem*{formula}{{Main Formula}}

\theoremstyle{definition}
\newtheorem*{remark}{Remark} 
 
\begin{document}

\baselineskip=17pt
\title[The first Chern class of the Verlinde bundles]{The first Chern class of the Verlinde bundles}
\author {A. Marian}
\address {Department of Mathematics, Northeastern University}
\email {a.marian@neu.edu}
\author {D. Oprea}
\address {Department of Mathematics, University of California, San Diego}
\email {doprea@math.ucsd.edu}
\author{R. Pandharipande}
\address{Department of Mathematics, ETH Z\"{u}rich}
\email {rahul@math.ethz.ch}
\dedicatory{To the memory of F. Hirzebruch}
\begin{abstract}
A formula for the first Chern class of the Verlinde bundle over the moduli space of smooth genus $g$ curves is given. A finite-dimensional argument is presented in rank $2$ using  geometric symmetries obtained from strange duality, relative Serre duality, and Wirtinger duality together with the projective flatness of the Hitchin connection. A derivation using conformal-block methods is presented in higher rank. An expression for the first Chern class over the compact moduli space of curves is obtained. 
\end {abstract}

\maketitle

\setcounter{tocdepth}{1}
\tableofcontents

\section {Introduction} 
\subsection {The slopes of the Verlinde complexes} Let $\mathcal M_g$ be the moduli space
of nonsingular curves of genus $g\geq 2$.
Over $\mathcal M_g$, we consider the relative moduli space of rank $r$ slope-semistable bundles 
of degree $r(g-1)$,
$$\nu: \mathcal {U}_g(r, r(g-1))\to \mathcal M_g\ .$$  The moduli space comes equipped with a canonical universal theta bundle corresponding to the divisorial locus $$\Theta_{r}=\{(C, E): h^0(E)\neq 0\}.$$ Pushing forward the pluritheta series, we obtain a canonical Verlinde complex \footnote{To avoid technical difficulties, it will be convenient to use the coarse moduli schemes of semistable vector bundles throughout most of the paper. Nonetheless, working over the moduli stack yields an equivalent definition of the Verlinde complexes, see Proposition $8.4$ of \cite{beauvillelaszlo}.}$$\mathbb V_{r, k}=\mathbf R\nu_{\star} \left(\Theta^k_{r}\right)$$ over $\mathcal M_g$. For $k\geq 1$, $\mathbb V_{r, k}$ is a vector bundle. 

The Verlinde bundles are known to be projectively flat \cite {Hi}. Therefore, their Chern characters satisfy the identity \begin{equation}\label{fl}\text{ch} (\mathbb V_{r, k})=\text{ rank } \mathbb V_{r, k} \ \cdot\ \exp \left( \frac{c_1(\mathbb V_{r, k})}{\text{rank } \mathbb V_{r, k}}\right).\end{equation}
The rank of ${\mathbb V}_{r,k}$ is given by the well-known Verlinde formula, see \cite {beauville}. We are interested here in calculating the slope $$\mu(\mathbb V_{r, k})=\frac{c_1(\mathbb V_{r, k})}{\text{rank } \mathbb V_{r, k}} 
\in H^2(\m_g,\mathbb{Q}).$$ Since the Picard rank of $\mathcal M_g$ is 1, we can express the slope in the form $$\mu(\mathbb V_{r, k})=s_{r, k}\, \lambda$$ 
where $\lambda \in H^2 (\m_g, {\mathbb Q})$ is the first Chern class of the Hodge bundle. We seek to determine the rational numbers $s_{r, k}\in \mathbb Q.$ 
By Grothendieck-Riemann-Roch for the push-forward defining the Verlinde bundle, $s_{r, k}$ is in fact a rational function in $k$. 

\begin {formula}\label{conj1}The Verlinde slope is \begin{equation}\label{mainformula}\mu(\mathbb V_{r, k})=\frac{r(k^2-1)}{2(k+r)}\, \lambda\ .\end{equation}
\end {formula} 
\noindent 

The volume of the moduli space $\mathcal {U}_C(r, r(g-1))$ 
of bundles over a fixed curve with respect to the symplectic form induced by the canonical theta divisor is known to be given in terms of the irreducible representations $\chi$ of the group $SU_r$ : $$\text{vol}_r=\int_{\u_C(r, r(g-1))} \exp(\Theta)=c_r\cdot \sum_{\chi} \left(\frac{1}{\dim \chi}\right)^{2g-2}$$ for the constant $$c_r=(2\pi)^{-r(r-1)(g-1)}(1!\, 2!\, \cdots\,  (r-1)!)^{-(g-1)}.$$ Taking the $k\to \infty$ asymptotics in formula \eqref{mainformula} and using \eqref{fl}, we obtain as a consequence an expression for the cohomological push-forward:
$$\nu_{\star} \left(\exp({\Theta})\right)=\text{vol}_r\cdot \exp \left(\frac{r}{2}\lambda\right)\ .$$ This is a higher rank generalization of an equality over the relative Jacobian observed in \cite {vdG}. 

\vskip.1in

In {\bf Part I} of this paper, we are concerned with a finite-dimensional geometric proof of the Main Formula. In {\bf Part II}, we give a derivation via conformal blocks. We also extend the formula over the boundary of the moduli space. Let us now detail the discussion.    
\vskip.1in

For the finite dimensional argument, we note four basic symmetries of the geometry:
\begin {itemize}
\item [(i)] Relative level-rank duality for the moduli space of bundles over $\mathcal M_g$ will be shown to give the identity $$s_{r,k} + s_{k,r} = \frac{kr-1}{2}\ .$$
\item [(ii)] Relative duality along the the fibers of $\mathcal {SU}_g(r, {\mathcal O})\to \mathcal M_g$ leads to $$s_{r,k} + s_{r, -k -2r} = -2r^2.$$
\item [(iii)] The initial conditions in rank $1$, and in level $0$ are $$\mu({\mathbb V}_{1,k}) = \frac{k-1}{2}\ ,\,\,\, \ \ \mu(\mathbb V_{r, 0})=-\frac{1}{2}
\ .$$
\item [(iv)] The projective flatness of the Verlinde bundle.
\end{itemize} The four features of the geometry will be shown to determine the Verlinde slopes completely in the rank 2 case, proving:

\begin {theorem} \label{main} The Verlinde bundle ${\mathbb V}_{2,k}$ has slope $$\mu(\mathbb V_{2, k})=\frac{k^2-1}{k+2}\, \lambda\ .$$
\end {theorem} 

In arbitrary rank, the symmetries entirely determine the slopes in the Main Formula \eqref{mainformula} under one additional assumption.  This assumption concerns the roots of the Verlinde polynomial $$v_g(k)=\chi(\s\u_C(r, \mathcal O), \Theta^{k})$$ giving the $SU_r$ Verlinde numbers at level $k$.  Specifically, with the exception of the root $k=-r$ which should have multiplicity exactly $(r-1)(g-1)$, all the other roots of $v_g(k)$ should have multiplicity less than $g-2$. Numerical evidence suggests this is true. 
\vskip.1in

Over a fixed curve $C$, the moduli spaces of bundles with fixed determinant $\s\u_C(2r, \mathcal O_C)$ and $\s\u_C(2r, \omega_C^r)$ are isomorphic. Relatively over $\mathcal M_g$ such an isomorphism does not hold. 
Letting $\Theta$ denote the canonical theta divisor in $$\nu: {\mathcal {SU}}_{g}(2r, \omega^r) \to \m_g\ ,$$ we may investigate the slope of $$\mathbb W_{2r, k}=\mathbf R\nu_{\star} (\Theta^k)\ .$$ The following statement is equivalent to Main Formula \eqref{mainformula} via Proposition \ref{slid2} of Section \ref{candet}. 
As will be clear in the proof, the equivalence of the two statements corresponds geometrically to the relative version of Wirtinger's duality for level $2$ theta functions. 
\begin {theorem} The Verlinde bundle $\mathbb W_{2r, k}$ has slope
$$\mu(\mathbb W_{2r, k})=\frac{k(2rk+1)}{2(k+2r)}\, \lambda\ .$$
\end {theorem} 
\noindent

\vskip.1in

In {\bf Part II}, we deduce the Main Formula from a representation-theoretic perspective by connecting results in the conformal-block literature. In particular, essential to the derivation are the main statements in \cite{T}. There, an action of a suitable Atiyah algebra, an analogue of a sheaf of differential operators, is used to describe the projectively flat WZW connection.  Next, results of Laszlo \cite{L} identify conformal blocks and the bundles of theta functions aside from a normalization ambiguity. An integrality argument fixes the variation over moduli of the results of \cite {L}, yielding the main slope formula. This is explained in Section \ref{cfbl}. 

Finally, in the last section, we consider the extension of the Verlinde bundle over the compact moduli space $\overline {\mathcal M}_g$ via conformal blocks. The Hitchin connection is known to acquire regular singularities along the boundary \cite {TUY}.  The formulas for the first Chern classes of the bundles of conformal blocks are given in Theorem \ref{form} of Section \ref{extension}. They specialize to the genus $0$ expressions of \cite {F} in the simplified form of \cite {Mu}. 
\vskip.1in
\noindent {\bf Related work.} In genus $0$, the conformal block bundles have been studied in recent years in connection to the nef cone of the moduli space $\overline{\m}_{0,n}$, see \cite {AGS}, \cite {AGSS}, \cite{F}, \cite {Fe}, \cite {GG}, \cite{Sw}. In higher genus, the conformal block bundles have been considered in \cite {S} in order to study certain representations arising from Lefschetz pencils. The method of \cite{S} is to use Segal's loop-group results. Unfortunately, the geometry underlying \cite{S} is not uniquely specified.

There are at least two perspectives on the study of the higher Chern classes of the
Verlinde bundle. Via a version of Thaddeus wall-crossing studied relatively over $\mathcal M_{g,1}$,
an approach to the higher Chern class of the Verlinde bundle is pursued in \cite {FMP}.
Projective flatness then yields nontrivial relations in the tautological ring $R^{\star}(\mathcal M_{g, 1})$ of
the moduli space of curves. Whether these relations always lie in the Faber-Zagier set \cite {PP} is an open question.

A completely different point of view is taken in \cite{MOPPZ}. The Chern character of the
conformal block bundle defines a semisimple CohFT via the fusion rules. The Givental- Teleman theory provides a classification up to an action of the Givental group. A unique element of the classification is selected by the projective flatness condition and the first Chern class calculation. The outcome is a clean formula for the higher Chern classes extending the first Chern class result of Theorem 3 proven here. However, since the latter formula incorporates the projective flatness as an input, no nontrivial relations in $R^{\star}(\mathcal M_{g,1})$ are obtained.

\subsection {Acknowlegements}
We thank Carel Faber for the related computations in 
\cite{FMP} and Ivan Smith for correspondence concerning \cite{S}.
Our research was furthered
during the {\em Conference on Algebraic Geometry} in
July 2013 at the University of Amsterdam. We thank the
organizers for the very pleasant environment.

A.M. and D.O. were partially supported by the NSF grants DMS 1001604, DMS 1001486, DMS 1150675, as well as by Sloan Foundation Fellowships. R.P. was partially supported by the grant ERC-2012-AdG-320368-MCSK.
\vskip.1in
\section*{Part I: Finite-dimensional methods}

\section{Jacobian geometry}

In this section, we record useful aspects of the geometry of relative Jacobians over the moduli space of curves. The results will be used to derive the slope identities of  Section \ref{slids}. 

Let $\m_{g,1}$ be the moduli space of nonsingular 1-pointed genus $g\geq 2$ curves, 
and let $$\pi: {\mathcal C} \to \m_{g,1}, \, \, \, \, \, \, \sigma: \m_{g,1} \to {\mathcal C}$$ be the universal curve and the tautological section respectively. We set $\gbar = g-1$ for convenience. The following line bundle will play an important role in subsequent calculations:
$${\mathcal L} \to \m_{g,1}, \, \, \,\ \ \  {\mathcal L} = \left (  \det {\bf R}\pi_{\star} {\mathcal O_{\mathcal C} (\gbar \sigma)} \right )^{-1}.$$
An elementary Grothendieck-Riemann-Roch computation applied to the
morphism 
$\pi$ yields  
$$c_1 ({\mathcal L}) = -\lambda + \binom{g}{2} {\Psi},$$ 
where $$\Psi\in H^2(\mathcal M_{g, 1},\mathbb{Q})$$ is the cotangent class. 

Consider $p: {\mathcal J} \to \m_{g,1}$ the relative Jacobian of degree 0 line bundles. We let $$\widehat\Theta \to {\mathcal J}$$ be the line bundle associated 
to the divisor 
\begin{equation}
\label{thetabasic}
\{ (C, p, L) \, \, \, \text{with} \, \, \, H^0 (C, L (\gbar\, p)) \neq 0 \},
\end{equation}
and let $$\theta = c_1 (\widehat \Theta )$$ be the corresponding
 divisor class. We show

\begin{lemma}\label{lem1}
$p_{\star} \left(e^{n\theta}\right) = n^g e^{\frac{n \lambda}{2}}$.
\end{lemma}
\proof Since the pushforward sheaf $p_{\star} (\widehat\Theta)$ has rank 1 and a nowhere-vanishing section obtained from the divisor \eqref{thetabasic}, we
see that $$p_{\star} \left(\widehat\Theta\right) = {\mathcal O}_{\m_{g,1}}.$$ 
The relative tangent bundle of $$p: {\mathcal J} \to \m_{g,1}$$ is the pullback of the dual Hodge bundle $\mathbb E^{\vee} \to \m_{g,1}$, with Todd genus $$\text{Todd} \, {\mathbb E}^{\vee} = e^{-\frac{\lambda}{2}},$$ 
see \cite {vdG}. Hence, Grothendieck-Riemann-Roch yields 
$$p_{\star} (e^{\theta}) = e^{\frac{\lambda}{2}}.$$
The Lemma  follows immediately. \qed

\vskip.1in

Via Grothendieck-Riemann-Roch for $p_{\star} \left(\widehat \Theta^k\right)$,
 we obtain the following corollary of Lemma \ref{lem1}.

\begin{corollary} We have $$s_{1, k}=\frac{k-1}{2}.$$ 
\end{corollary}

\vskip.1in

We will later require the following result obtained as a consequence
of Wirtinger duality. Let $(-1)^\star\theta$ denote the pull-back of $\theta$ by the
involution $-1$ in the fibers of $p$.

\begin{lemma}\label{slcomp} $p_{\star} \left(e^{n(\theta+(-1)^{\star}\theta)}\right)=(2n)^{g} e^{2n\,
 c_1(\mathcal L)}.$
\end{lemma}
\proof We begin by recalling the classical Wirtinger duality for level $2$ theta functions. For a principally polarized abelian variety $(A, \widehat \Theta)$, we consider the map $$\mu:A\times A\to A\times A$$ given by $$\mu(a, b)=(a+b, a-b).$$ We calculate the pullback line bundle \begin{equation}\label{identity}\mu^{\star}(\widehat \Theta\boxtimes \widehat \Theta)=\widehat \Theta^2\boxtimes (\widehat \Theta\otimes (-1)^{\star}\widehat \Theta).\end{equation} The unique section of $\widehat \Theta\boxtimes \widehat \Theta$ gives a natural section of the bundle \eqref{identity}, inducing by K\"unneth decomposition an isomorphism $$H^0(A, \widehat \Theta^2)^{\vee}\to H^0(A, \widehat \Theta\otimes (-1)^{\star}\widehat \Theta),$$ see \cite {M}. 

We carry out the same construction for the relative Jacobian $$\mathcal J\to \mathcal M_{g, 1}.$$ Concretely, we let $$\mu:\mathcal J\times_{\mathcal M_{g, 1}} \mathcal J\to \mathcal J\times_{\mathcal M_{g, 1}} \mathcal J$$ be relative version of the map above.
The fiberwise identity \eqref{identity} needs to be corrected by a line bundle twist from $\mathcal M_{g,1}$: \begin{equation}\label{wirt}\mu^{\star}(\widehat \Theta\boxtimes \widehat \Theta)=\widehat \Theta^2\boxtimes \left(\widehat \Theta\otimes (-1)^{\star}\widehat \Theta\right)\otimes \mathcal T.\end{equation} We determine $$\mathcal T=\mathcal L^{-2}$$ by constructing a section $$s:\mathcal M_{g, 1}\to \mathcal J\times_{\mathcal M_{g, 1}}\mathcal J,$$ for instance $$s(C, p)=(\mathcal O_C, \mathcal O_C).$$ Pullback of \eqref{wirt} by $s$ then gives the identity $$\mathcal L^2=\mathcal L^2\otimes \mathcal L^2\otimes \mathcal T$$ yielding the expression for $\mathcal T$ claimed above. 
Pushing forward \eqref{wirt} to $\mathcal M_{g, 1}$ we obtain the relative Wirtinger isomorphism $$\left(p_{\star} (\widehat \Theta^{2})\right)^{\vee}\cong p_{\star} \left(\widehat \Theta\otimes (-1)^{\star} \widehat \Theta\right)\otimes \mathcal L^{-2}.$$ 

We calculate the Chern characters of both bundles via Grothendieck-Riemann-Roch. We find $$\left(p_{\star}(e^{2\theta})e^{-\frac{\lambda}{2}}\right)^{\vee}=p_{\star} (e^{\theta+(-1)^{\star}\theta}) \cdot e^{-\frac{\lambda}{2}}\cdot e^{-2c_1(\mathcal L)}.$$ We have already seen that $$p_{\star}(e^{2\theta})=2^{g} e^{\lambda},$$ hence the above identity becomes $$p_{\star} (e^{\theta+(-1)^{\star} \theta})=2^{g}e^{2c_1(\mathcal L)}.$$ The formula in the Lemma follows immediately. \qed

\section{Slope identities} \label{slids}
\subsection{Notation}
In the course of the argument, we will occasionally view the spaces of bundles over the moduli space $\m_{g,1}$ of pointed genus $g$ curves:
$$\s\u_{g,1} (r, \mathcal O) = \s\u_g ( r, \mathcal O) \times_{\m_g} \m_{g,1}, \, \, \, \, \u_{g,1} (r, r \gbar) = \u_g (r, r \gbar) \times_{\m_g} \m_{g,1}.$$ To keep the notation simple, we will use $\nu$ to denote all bundle-forgetting maps from the relative moduli spaces of bundles to the space of (possibly pointed) 
nonsingular curves. 

Over the relative moduli space $\, \u_{g,1} (r, r \gbar)$ there is a natural determinant line bundle $$\Theta_r \to \u_{g,1} (r, r\gbar),$$ endowed with a canonical section vanishing on the locus $$\theta_r = \{E \to C \, \, \text{with} \, \, H^0 (C, \, E) \neq 0 \}.$$ 

We construct analogous theta bundles for the moduli space of bundles with trivial determinant, and decorate them with the superscript ``+" for clarity. Specifically, we consider the determinant line bundle and corresponding divisor 
$$\Theta_r^+ \to \s\u_{g,1} (r, \mathcal O), \, \, \, \, \theta_r^+=\{ (C, \, p,\,  E\to C)\, \, \text{with} \, \, H^0 (C, \, E (\gbar p) ) \neq 0\}.$$ Pushforward yields an associated Verlinde bundle $$\mathbb V_{r, k}^+=\mathbf R{\nu}_{\star} \left(\left(\Theta_{r}^+\right)^{k}\right)\to \mathcal M_{g, 1}.$$ This bundle is however not defined over the unpointed moduli space $\mathcal M_{g}$. 

While the first Chern class of $\mathbb V_{r, k}$ is necessarily a multiple of $\lambda$, the first Chern class of $\mathbb V_{r, k}^+$ is a combination of $\lambda$ and the cotangent class $$\Psi\in H^2(\mathcal M_{g, 1},\mathbb{Q}).$$ 

\subsection{Strange duality} Using a relative version of the level-rank duality over moduli spaces of bundles on a smooth curve, we first prove the following slope symmetry.
\begin{proposition} For any positive integers $k$ and $r$, we have
\label{verlindesym}
$$s_{k, r} + s_{r, k} = \frac{kr - 1}{2}.$$
\end{proposition}

\noindent {\it Proof.} Let $$\tau: \s\u_{g,1} (r, \mathcal O) \, \times_{\m_{g,1}} \, \u_{g,1} (k, k\gbar) \, \longrightarrow \,  \u_{g,1} (kr, kr\gbar)$$ be the tensor product map, $$\tau (E, F) = E \otimes F.$$ Over each fixed pointed curve $(C, p) \in \m_{g,1}$ we have, as explained for instance in \cite{beauville},  
\begin{equation}
\label{sd0}
\tau^{\star} \Theta_{kr} \simeq \left(\Theta_r^{+}\right)^{k} \boxtimes \Theta_k^r \, \, \, \, \, \, \, \, \text{on} \, \, \, \s\u_C (r, \mathcal O) \times \u_C (k, k\gbar).
\end{equation}
The natural divisor $$\tau^{\star} \theta_{kr} = \{(E, \, F) \, \, \text{with} \, \, H^0 (E \otimes F) \neq 0 \}$$ induces the strange duality map, defined up to multiplication by scalars,
\begin{equation}
\label{sd1}
H^0 \left  ( \s\u_C (r, \mathcal O), (\Theta_r^+)^k \right )^{\vee} \, \longrightarrow \, H^0 \left ( \u_C (k, k \gbar), \Theta_k^r \right ).
\end{equation}
This map is known to be an isomorphism \cite{Bel}, \cite{mo}, \cite {P}.

Relatively over $\m_{g,1}$ we write, using the fixed-curve pullback identity \eqref{sd0}, 
\begin{equation}
\label{sdrel0}
\tau^{\star} \Theta_{kr} \simeq \left(\Theta_r^{+}\right)^{k} \boxtimes \Theta_k^r \otimes \nu^{\star} {\mathcal T} \, \, \, \, \, \, \, \, \, \text{on} \, \, \, \s\u_{g,1} (r, \mathcal O) \, \times_{\m_{g,1}} \, \u_{g,1} (k, k\gbar),
\end{equation}
for a line bundle twist $${\mathcal T} \to \m_{g,1}.$$ We will determine $$\mathcal T = {\mathcal L}^{kr}, \, \, \, \text{so that} \, \, c_1(\mathcal T)=kr\left(\lambda - \binom{g}{2} {\Psi}\right).$$ 
To show this,  we pull back \eqref{sdrel0} via the section $$s: \m_{g,1} \to \s\u_{g,1} (r, \mathcal O) \, \times_{\m_{g,1}} \, \u_{g,1} (k, k\gbar), \, \, \, s (C, p) = ({\mathcal O}_C^{\oplus r}, \, \, {\mathcal O}_C (\gbar p)^{\oplus k}),$$
obtaining
$${\mathcal L}^{kr} \simeq {\mathcal L}^{kr} \otimes {\mathcal L}^{kr} \otimes {\mathcal T},$$ hence the claimed expression for $\mathcal T$.

Pushing forward \eqref{sdrel0} now, we note, as a consequence of \eqref{sd1}, the isomorphism of Verlinde vector bundles over $\m_{g,1}$,
$$\left({{\mathbb {V}}_{r,k}^+}\right)^{\vee} \simeq \mathbb V_{k,r} \otimes {\mathcal T}.$$
We conclude  
$$ - \mu \left ( {{\mathbb {V}}_{r,k}^+} \right ) = \mu \left ( {\mathbb V}_{k, r} \right ) + c_1 ({\mathcal T}),$$
hence 
$$ - \mu \left ( {{\mathbb {V}}_{r,k}^+} \right ) = \mu \left ( {\mathbb V}_{k, r} \right ) + kr \left ( \lambda - \binom{g}{2} {\Psi} \right ).$$
The equation, alongside the following Lemma, allows us to conclude Proposition \ref{verlindesym}. \qed 

\begin{lemma} We have
\label{slopecomparison}
\begin{eqnarray*}\mu \left ( {{\mathbb {V}}_{r,k}} \right ) &=& \mu \left ( {{\mathbb {V}}_{r,k}^+} \right ) + \frac{kr-1}{2} \lambda - kr c_1 ({\mathcal L})\\&=&\mu\left ( {{\mathbb {V}}_{r,k}^+} \right ) +\frac{3kr-1}{2}\lambda-kr \binom{g}{2}\Psi.\end{eqnarray*}
\end{lemma}

\noindent 
{\it Proof.}  To relate $ \mu \left ( {{\mathbb {V}}_{r,k}^+} \right )$ and $ \mu \left ( {{\mathbb {V}}_{r,k}} \right )$ we use a slightly twisted version of the tensor product map $\tau$ in the case $k=1$. More precisely we have the following diagram, where the top part is a fiber square
\begin{equation*}
\xymatrix{\s\u_{g,1} (r, \mathcal O) \times_{\m_{g,1}} {\mathcal J} \ar[r]^{\,\,\,\,\,\,\,\,\,\,\,t} \ar[d]^{\bar{q}} &  \u_g (r, r \gbar) \ar[d]^{q} \\ {\mathcal J} \ar[r]^{r} \ar[dr]^{p} & {\mathcal J} \ar[d]^{p} \\  & \m_{g,1} } . \end{equation*}
Here, as in the previous section, we write $$p: {\mathcal J} \to \m_{g,1}$$ for the relative Jacobian of degree 0 line bundles, while $r$ denotes multiplication by $r$ on ${\mathcal J}$. Furthermore, for a pointed curve $(C,p),$
$$t (E, L) = E \otimes L(\gbar p), \, \, \, \, q (E) = \det (E (-\gbar p)).$$ Finally, $\bar{q}$ is the projection onto ${\mathcal J}$.
The pullback equation \eqref{sdrel0} now reads
$$ t^{\star} \Theta_r \simeq \Theta_r^+ \boxtimes \widehat \Theta^r \otimes {\mathcal L}^{-r},$$
where, keeping with the previous notation, $\widehat \Theta \to {\mathcal J}$ is the theta line bundle associated with the divisor $$\theta : = \{(C, p, L\to C) \, \, \text{with} \, \, H^0 (C, L(\gbar p)) \neq 0 \}.$$
Using the pullback identity and the Cartesian diagram, we conclude
\begin{equation}
\label{sdrel00}
r^{\star} q_{\star} (\Theta_r^k) = \bar{q}_{\star} \left (\left ( \Theta_r^{+} \right )^k  \boxtimes \widehat \Theta^{kr} \otimes {\mathcal L}^{-kr} \right ) = p^{\star} {\mathbb V}_{r,k}^+ \otimes \widehat \Theta^{kr} \otimes {\mathcal L}^{-kr} \, \, \text{on} \, \, {\mathcal J}.
\end{equation}
We are however interested in calculating $$\text{ch} \, {\mathbb V}_{r,k} = \text{ch} \,\nu_{\star} \Theta_{r}^k = \text{ch} \, p_{\star}  ( q_{\star} \Theta_r^k  ).$$
We have recorded in Lemma \ref{lem1} the Todd genus of the 
the relative tangent bundle of $$p: {\mathcal J} \to \m_{g,1}$$ to be $$\text{Todd} \, {\mathbb E^{\vee}} = e^{-\frac{\lambda}{2}}.$$ 
Grothendieck-Riemann-Roch then gives 
$$\text{ch} \, {\mathbb V}_{r,k} = e^{-\frac{\lambda}{2}} \, p_{\star} ( \text{ch} \,( q_{\star} \Theta_r^k)).$$
We further write, on ${\mathcal J},$
$$\text{ch} \,( q_{\star} \Theta_r^k) = \frac{1}{r^{2g}}\,  r^{\star}  \text{ch} \,( q_{\star} \Theta_r^k) = \frac{1}{r^{2g}} \,  \text{ch} \,( r^{\star} q_{\star} \Theta_r^k) =  \frac{1}{r^{2g}} \, e^{kr \theta - kr c_1({\mathcal L})} \, p^{\star} \text{ch} \,  {\mathbb V}_{r,k}^+,$$
where \eqref{sdrel00} was used.
We obtain
$$\text{ch} \, {\mathbb V}_{r,k} = \frac{1}{r^{2g}}\, e^{-\frac{\lambda}{2}} \, e^{-kr c_1({\mathcal L})} \, \left ( p_{\star} e^{kr \theta} \right )  \text{ch} \,  {\mathbb V}_{r,k}^+ \, \, \, \text{on} \, \, \, \m_{g,1}.$$
The final $p$-pushforward in the identity above was calculated in Lemma \ref{lem1}. Substituting, we obtain
$$\text{ch} \, {\mathbb V}_{r,k} = \frac{k^g}{r^{g}}\, e^{\frac{(kr - 1) \lambda}{2}} \, e^{-kr c_1({\mathcal L})} \,   \text{ch} \,  {\mathbb V}_{r,k}^+ \, \, \, \text{on} \, \, \, \m_{g,1}.$$
Therefore, 
$$\mu \left (  {\mathbb V}_{r,k} \right ) = \mu \left ( {\mathbb V}_{r,k}^+ \right ) + \frac{kr-1}{2} \lambda - kr c_1 ({\mathcal L}),$$ which is the assertion of Lemma \ref{slopecomparison}.\qed

\subsection {Relative Serre duality}

We will presently deduce another identity satisfied by the numbers $s_{r,k}$ using relative Serre duality for the forgetful morphism $$\nu: \s\u_{g,1}(r, \mathcal O) \to \m_{g,1}\ .$$ 
\begin {proposition} \label{rs}We have $$s_{r, k}+s_{r, -k-2r}=-2r^2\ .$$
\end {proposition}

\proof By relative duality, we have $${\mathbb V}^+_{r, k}=\mathbf R\nu_{\star} \left(\left(\Theta_{r}^{+}\right)^{k}\right)\cong \mathbf R\nu_{\star} \left(\left(\Theta_{r}^{+}\right)^{-k}\otimes \mathcal \omega_{\nu}\right)^{\vee}[(r^2-1)(g-1)]\ .$$ 
We determine the relative dualizing sheaf of the morphism $\nu$. As explained in Theorem E of \cite {DN}, the fibers of the morphism $$\nu:\mathcal {SU}_{g,1}(r, \mathcal O)\to \mathcal M_{g,1}$$ are Gorenstein, hence the relative dualizing sheaf is a line bundle. Furthermore, along the fibers of $\nu$, the canonical bundle equals $\left ( \Theta_r^{+} \right )^{-2r}.$ Thus, up to a line bundle twist $\mathcal T\to \mathcal M_{g,1}$, we have 
\begin{equation}\label{twisto}\omega_{\nu}=\left (\Theta_r^{+} \right )^{-2r}\otimes \nu^{\star} \mathcal T.\end{equation}The twist $\mathcal T$ will be found via a Chern class calculation to be
$$c_1(\mathcal T)= -(r^2+1)\lambda+ 2r^2 \binom{g}{2} {\Psi}.$$ 
Since $$\left({\mathbb V}_{r, k}^{+}\right)^{\vee}\cong \mathbf R\nu_{\star} ((\Theta_{r}^{+})^{-k}\otimes \mathcal \omega_{\nu})[(r^2-1)(g-1)]=\mathbb V_{r, -k-2r}^+\otimes \mathcal T[(r^2-1)(g-1)],$$ we obtain taking slopes that $$-\mu(\mathbb V_{r, k}^+)=\mu(\mathbb V_{r, -k-2r}^+)+\left(-(r^2+1)\lambda+ 2r^2 \binom{g}{2} {\Psi}\right).$$ The proof is concluded using Lemma \ref{slopecomparison}.

To determine the twist $\mathcal T$, we begin by restricting \eqref{twisto} to the smooth stable locus of the moduli space of bundles $$\nu:\mathcal {SU}^{s}_{g,1}(r,\mathcal O)\to \mathcal M_{g, 1}.$$ There, the relative dualizing sheaf is the dual determinant of the relative tangent bundle. By Corollary $4.3$ of \cite {DN}, adapted to the relative situation, the Picard group of the coarse moduli space and the Picard group of the moduli stack are naturally isomorphic. We therefore consider \eqref{twisto} over the moduli stack of stable bundles. (We do not introduce separate notation for the stack, for simplicity.) Let $$\mathcal E\to \mathcal {SU}^{s}_{g,1}(r,\mathcal O)\times_{\mathcal M_{g, 1}}\mathcal C$$ denote the universal vector bundle of rank $r$ over the stable part of the moduli stack. We write $$\pi:\mathcal {SU}^{s}_{g,1}(r, \mathcal O)\times_{\mathcal M_{g, 1}}\mathcal C\to \mathcal {SU}^s_{g,1}(r, \mathcal O)$$ for the natural projection. Clearly, $$\Theta_r^{+}=(\det \mathbf R \pi_{\star} (\mathcal E \otimes {\mathcal O}_{\mathcal C} (\gbar \sigma)))^{-1}.$$ The relative dualizing sheaf of the morphism $\nu$ is expressed as $$\omega_{\nu}=\mathbf R \pi_{\star} {Hom}(\mathcal E, \mathcal E)_{(0)}=\mathbf R\pi_{\star}{Hom}(\mathcal E, \mathcal E)-\mathbf R \pi_{\star} \mathcal O.$$ We therefore have $$c_1(\omega_{\nu})=c_1(\mathbf R\pi_{\star}{Hom}(\mathcal E, \mathcal E))-\lambda.$$ Using $\omega$ for the relative dualizing sheaf along the fibers of $\pi$, we calculate \begin{eqnarray*}c_1(\omega_{\nu})&+&2rc_1(\Theta_r^{+})=c_1\left(\mathbf R\pi_{\star} \it{End } \mathcal E\right)-\lambda-2rc_1\left(\mathbf R\pi_{\star} (\mathcal E \otimes {\mathcal O}_{\mathcal C} (\gbar \sigma)) \right)\\&=&\pi_{\star} \left[\left(1-\frac{\omega}{2}+\frac{\omega^2}{12} \right) \left(r^2+((r-1)c_1(\mathcal E)^2-2rc_2(\mathcal E)\right)\right]_{(2)}-\lambda\\&-&2r\pi_{\star}\left[\left(1-\frac{\omega}{2}+\frac{\omega^2}{12}\right)\left(r+c_1(\mathcal E)+\frac{1}{2}c_1(\mathcal E)^2-c_2(\mathcal E)\right) (1 + \gbar \sigma - \frac{\gbar^2}{2} \sigma \, {\Psi} )\right]_{(2)}\\&=&-(r^2+1)\lambda+ 2r^2 \binom{g}{2} {\Psi} + \pi_{\star}(r\omega\cdot c_1(\mathcal E)- 2r \gbar \,\sigma \cdot c_1 (\mathcal E) - c_1(\mathcal E)^2).\end{eqnarray*} Since the determinant of $\mathcal E$ is trivial on the fibers of $\pi$, we may write $$\det \mathcal E=\pi^{\star}\mathcal A$$ for a line bundle $\mathcal A \to \mathcal {SU}^{s}_{g,1}(r, \mathcal O)$ with first Chern class
$$\alpha = c_1(\mathcal A).$$ We calculate $$\pi_{\star}(r\omega\cdot c_1(\mathcal E)- 2r \gbar \,\sigma \cdot c_1 (\mathcal E) - c_1(\mathcal E)^2) = 2r \gbar \alpha - 2r \gbar \alpha - \pi_{\star} (\alpha^2) = 0,$$ and conclude $$\nu^{\star} c_1(\mathcal T)=c_1(\omega_{\nu})+2rc_1({\Theta_r^{+}}) = -(r^2+1)\lambda+ 2r^2 \binom{g}{2} {\Psi}.$$ This equality holds in the Picard group of the stable locus of the moduli stack and of the coarse moduli space. Since the strictly semistables have codimension at least $2$, the equality extends to the entire coarse space $\s\u_{g,1}(r, \mathcal O)$. Finally, pushing forward to $\mathcal M_{g, 1}$, we find the expression for the twist $\mathcal T$ claimed above.  
\qed

\subsection{Initial conditions} The next calculation plays a basic role in our argument. 
 
 \begin {lemma} \label{ic} We have $$s_{r, 0}=-\frac{1}{2} .$$
 \end{lemma}
 \proof Since the Verlinde number for $k=0$ over the moduli space $\u_C(r, r\gbar)$ is zero, the slope appears to have poles if computed directly. Instead,
we  carry out the calculation via the fixed determinant moduli space. 
The trivial bundle has no higher cohomology along the fibers of $$\nu:\mathcal {SU}_{g, 1} (\,r, \mathcal O)\to \mathcal M_{g,1}$$ by Kodaira vanishing. To apply the vanishing theorem, we use that the fibers of $\nu$ have rational singularities, and the expression of the dualizing sheaf of Proposition \ref{rs}.
Hence, $$\nu_{\star}\left(\mathcal O\right)=\mathcal O_{\mathcal M_{g,1}}.$$ Therefore $$\mu(\mathbb V_{r, 0}^+)=0$$ which then immediately implies 
$s_{r, 0}=-\frac{1}{2}$ by Lemma \ref{slopecomparison}.
 \qed
 
\subsection{Pluricanonical determinant} 
\label{candet}
We have already investigated
moduli spaces of bundles with trivial determinant. 
Here, we assume that the determinant is
  of degree equal to rank times $\overline{g}$ 
and
is a multiple of the canonical bundle.
The conditions 
require the rank to be even. 
Thus, we are concerned with the slopes of the complexes $$\mathbb W_{2r, k}={\mathbf R}\nu_{\star} \left(\Theta_{2r}^{k}\right),$$ where $$\nu:\s\u_{g}(2r, \omega^r)\to \mathcal M_g.$$ The following slope identity is similar to that of Lemma \ref{slopecomparison}:
\begin {proposition}\label{slid2}
We have
 $$\mu(\mathbb W_{2r, k})=\mu(\mathbb V_{2r,k})+\frac{\lambda}{2}\ .$$ In particular, via Theorem \ref{main}, we have $$\mu(\mathbb W_{2, k})=\frac{k(2k+1)}{2(k+2)}\lambda\ .$$
\end{proposition} 

\proof Just as in the proof of Lemma \ref{slopecomparison}, we relate $ \mu \left ( {{\mathbb {W}}_{2r,k}} \right )$ and $ \mu \left ( {{\mathbb {V}}_{2r,k}} \right )$ via the tensor product map $t$:
\begin{equation*}
\xymatrix{\s\u_{g,1} (2r, \omega^r) \times_{\m_{g,1}} {\mathcal J} \ar[r]^{\,\,\,\,\,\,\,\,\,\,\,\,\,\,\,\,t} \ar[d]^{\bar{q}} &  \u_g (2r, 2r \gbar) \ar[d]^{q} \\ {\mathcal J} \ar[r]^{2r} \ar[dr]^{p} & {\mathcal J} \ar[d]^{p} \\  & \m_{g,1} } . \end{equation*} We keep the same notation as in  Lemma  \ref{slopecomparison},
letting $$p: {\mathcal J} \to \m_{g,1}$$ denote the relative Jacobian of degree 0 line bundles, and writing $2r$ for the multiplication by $2r$ on ${\mathcal J}$. Furthermore, for a pointed curve $(C,p),$
$$t (E, L) = E \otimes L, \, \, \, \, q (E) = \det E \otimes \omega_C^{-r}$$ Finally, $\bar{q}$ is the projection onto ${\mathcal J}$. Recall that $\widehat \Theta$ denotes the theta line bundle on the relative Jacobian associated with the divisor $$\theta : = \{(C, p, L) \, \, \text{with} \, \, H^0 (C, L(\gbar p)) \neq 0 \}.$$ It is clear that $(-1)^{\star}\widehat \Theta$ has the associated divisor $$(-1)^{\star} \theta=\{(C, p, L) \, \, \text{with} \, \, H^0 (C, L\otimes \omega_C(-\gbar p)) \neq 0\}.$$

For a fixed pointed curve $(C, p)$, we have the fiberwise identity 
$$t^{\star}\Theta_{2r}=\Theta_{2r}\boxtimes \left(\widehat \Theta\otimes (-1)^{\star}\widehat \Theta\right)^r$$ on $\s\u_C(2r, \omega^r)\times \mathcal {J}_C$. Relatively over $\mathcal M_{g, 1}$, the same equation holds true up to a twist $\mathcal T\to \mathcal M_{g,1}$:
$$ t^{\star} \Theta_{2r} \simeq \Theta_{2r} \boxtimes \left(\widehat \Theta\otimes (-1)^{\star}\widehat \Theta\right)^r \otimes {\mathcal T}.$$ We claim that $$\mathcal T=\mathcal L^{-2r}.$$ Indeed, the twist can be found in the usual way, using a suitable section $$s:\mathcal M_{g, 1}\to \s\u_{g,1} (2r, \omega^r) \times_{\m_{g,1}} {\mathcal J},$$ for instance $$s(C, p)=(\omega_C^r(-\gbar p)^{\oplus r}\oplus \mathcal O_C(\gbar p)^{\oplus r}, \mathcal O_C).$$ Pulling back by $s$, we obtain the identity $$\left (\mathcal L\otimes \mathcal M \right)^r =  \left (\mathcal L\otimes \mathcal M \right )^{r} \otimes \left (\mathcal L\otimes \mathcal M \right )^{r}\otimes \mathcal T$$ where $$\mathcal L=\det \left(\mathbf R\pi_{\star} (\mathcal O_{\mathcal C}(\gbar \sigma))\right)^{-1},\,\, \mathcal M=\det \left(\mathbf R\pi_{\star} \,(\mathcal \omega_{\mathcal C}(-\gbar \sigma))\right)^{-1}.$$ In fact, by relative duality, $\mathcal M\cong \mathcal L$, so we conclude $$\mathcal T=\mathcal L^{-2r}.$$

Using the pullback identity and the Cartesian diagram, we find that over $\mathcal J$ we have
\begin{eqnarray}\label{plbck}
(2r)^{\star} q_{\star} \Theta_{2r}^k &=& \bar{q}_{\star} \left (\Theta_{2r}^k  \boxtimes \left(\widehat \Theta\otimes (-1)^{\star}\widehat \Theta \right)^{kr} \otimes {\mathcal L}^{-2kr} \right )\\ &=&\nonumber p^{\star} {\mathbb W}_{2r,k} \otimes \left(\widehat \Theta\otimes (-1)^{\star}\widehat \Theta\right)^{kr} \otimes {\mathcal L}^{-2kr}
\end{eqnarray}
Next, we calculate $$\text{ch} \, {\mathbb V}_{2r,k} = \text{ch} \,\nu_{\star} \Theta_{2r}^k = \text{ch} \, p_{\star}  ( q_{\star} \Theta_{2r}^k  )$$ via Grothendieck-Riemann-Roch:
$$\text{ch} \, {\mathbb V}_{2r,k} = e^{-\frac{\lambda}{2}} \, p_{\star} ( \text{ch} \,( q_{\star} \Theta_{2r}^k)).$$
We further evaluate, on ${\mathcal J},$
$$\text{ch} \,( q_{\star} \Theta_{2r}^k) = \frac{1}{(2r)^{2g}}\,  (2r)^{\star}  \text{ch} \,( q_{\star} \Theta_{2r}^k) = \frac{1}{(2r)^{2g}} \,  \text{ch} \,( (2r)^{\star} q_{\star} \Theta_r^k) $$ $$=  \frac{1}{(2r)^{2g}} \, e^{kr (\theta+(-1)^{\star}\theta) - 2kr c_1({\mathcal L})} \, p^{\star} \text{ch} \,  {\mathbb W}_{2r,k},$$
where \eqref{plbck} was used.
We obtain
$$\text{ch} \, {\mathbb V}_{2r,k} = \frac{1}{(2r)^{2g}}\, e^{-\frac{\lambda}{2}} \, e^{-2kr c_1({\mathcal L})} \, \left ( p_{\star} e^{kr (\theta+(-1)^{\star}\theta)} \right )  \text{ch} \,  {\mathbb W}_{2r,k} \, \, \, \text{on} \, \, \, \m_{g,1}.$$
The $p$-pushforward in the identity above is given by Lemma \ref{slcomp}. Substituting, we find $$\text{ch} \, {\mathbb V}_{2r,k} =\left(\frac{k}{2r}\right)^{g}\, e^{-\frac{\lambda}{2}}  \text{ch} \,  {\mathbb W}_{2r,k},$$ and taking slopes it follows that 
$$\mu(\mathbb V_{2r, k})=\mu(\mathbb W_{2r, k})-\frac{\lambda}{2}\ .$$\qed

\section{Projective flatness and the rank two case} 
\subsection{Projective flatness}\label{pf}
By the Grothendieck-Riemann-Roch theorem for singular varieties due to Baum-Fulton-MacPherson \cite {BFM}, the Chern character of $\mathbb V_{r, k}$ is a polynomial in $k$ with entries in the cohomology classes of $\mathcal M_g$. (Alternatively, we may transfer the calculation to a smooth moduli space of degree $1$ bundles using a Hecke modification at a point as in \cite {BS}, and then invoke the usual Grothendieck-Riemann-Roch theorem.) Taking account of the projective flatness identity \eqref{fl},
$$\text{ch} (\mathbb V_{r, k})=\text{rank }\mathbb V_{r, k} \ \cdot \
\exp \left(s_{r, k} \lambda\right),$$
we therefore write
$$\text{ch}_i ({\mathbb V}_{r,k}) = \sum_{j=0}^{r^2 \gbar + i+1} k^j \alpha_{i,j} = 
(\text{rank } \mathbb V_{r, k})\,  \frac{s_{r, k}^{i}}{i!}\, \lambda^{i} \, \, \,  \text{for} \, \, \, \ i\geq 0, \, \,\ \alpha_{i,j} \in H^{2i} (\m_g).$$ As the Vandermonde determinant is nonzero, for each $i$ we can express $\alpha_{i,j}$ in terms of $\lambda^i$. Since $\lambda^{g-2} \neq 0,$ we deduce that 
$$(\text{rank } \mathbb V_{r, k})\, s_{r, k}^{i} \ , \, \, \, \, \, 0 \leq i \leq g-2,$$ is a polynomial in $k$ of degree $r^2 \gbar + i+1,$ with coefficients that may depend on $r$ and $g$. The following is now immediate:
\begin{itemize}
\item [(i)] 
For each $r$ we can write $$s_{r,k} = \frac{a_r(k)}{b_r(k)}$$ as quotient of polynomials of minimal degree, with $$\deg a_r (k) -\deg b_r (k) \leq 1.$$
Setting $v_{g,r} (k) =  \text{rank } \mathbb V_{r, k},$ we also have $$b_r(k)^{g-2} \text { divides } v_{g,r} (k)$$ as polynomials in $\mathbb Q[k]$.
\end{itemize}
In addition, the following properties of the function $s_{r, k}$ have been established in the previous sections:
\begin {itemize}
\item [(ii)] $s_{1, k}=\frac{k-1}{2}$, $s_{r, 0}=-\frac{1}{2}$,
\item [(iii)] $s_{r, k}+s_{k, r}=\frac{kr-1}{2}$ for all $k, r\geq 1$,
\item [(iv)] $s_{r, k}+s_{r, -k-2r}=-2r^2$ for all $r\geq 1$ and all $k$. 
\end {itemize} 

\vskip.1in

Clearly, the function $$s_{r, k}=\frac{r(k^2-1)}{2(k+r)}$$ of formula \eqref{mainformula} satisfies symmetries (ii)-(iv). Therefore, the shift $$s'_{r, k}=s_{r, k}-\frac{r(k^2-1)}{2(k+r)}$$ satisfies properties similar to (i)-(iv):
\begin {itemize}
\item [(i)$'$] $s'_{r, k}$ is a rational function of $k$,
\item [(ii)$'$] $s'_{1, k}=0$ for all $k$, and $s'_{r, 0}=0$ for all $r\geq 1$,
\item [(iii)$'$] $s'_{r, k}+s'_{k, r}=0$ for $r, k\geq 1$,
\item [(iv)$'$] $s'_{r, k}+s'_{r, -k-2r}=0$ for all $r\geq 1$ and all $k$.
\end {itemize}
\vskip.1in
\subsection{The rank two analysis}
To prove Theorem \ref{main},  we now show that $s'_{2, k}=0$ for all $k.$ Of course, $s'_{2, 0}=0$ by (ii)$'$. Also by (ii)$'$, we know that $s'_{1, 2}=0$, hence by (iii)$'$ we find $$s'_{2, 1}=0.$$ Similarly, $$s'_{2, 2}=0$$ also by (iii)$'$. 
Using (iv)$'$, we obtain that $$s'_{2, 0}=s'_{2, 1}=s'_{2, 2}=s'_{2, -4}=s'_{2, -5}=s'_{2, -6}=0.$$ Finally, we make use of the projective flatness of $\mathbb V_{2, k}.$ 
The Verlinde formula reads \cite {beauville} $$v_{g,2} (k)=k^{g}\left(\frac{k+2}{2}\right)^{g-1} \left(\sum_{j=1}^{k+1} \frac{1}{\sin^{2g-2}\frac{j\pi}{k+2}}\right).$$ The polynomial $v_{g,2} (k)$ admits $k=0$ as a root of order $g$ and $k=-2$ as a root of order $(g-1)$. Indeed, it was shown by Zagier that 
$$\widehat v_{g}(k+2)=\sum_{j=1}^{k+1} \left(\frac{1}{\sin \frac{j\pi}{k+2}}\right)^{2g-2}$$ is a polynomial in $k+2$ such that 
$$\widehat v_{g}(0)<0,$$ see Remark 1 on page 4 of \cite{Z}.

Let us write $$b_2(k)=(k+2)^{m}k^n B(k)$$ for a polynomial $B$ which does not have $0$ and $-2$ as roots. By property (i) above, we obtain $$m\leq \frac{g-1}{g-2}\implies m\leq 1.$$ Similarly $$n\leq \frac{g}{g-2}\implies n\leq1$$ unless $g=3,4$. Also, $B(k)^{g-2}$ divides the Verlinde polynomial $\widehat v_{g}(k+2)$ which has degree $4g-3-(g-1)-g=2g-2$. Thus $$(g-2)\deg B\leq 2g-2\implies \deg B\leq 2$$ except possibly when $g=3, 4$. 
In conclusion $$s'_{2, k}=\frac{a_2(k)}{B(k)(k+2)^mk^n}-\frac{k^2-1}{k+2}=\frac{A(k)}{B(k)(k+2)k}$$ for  a polynomial $$A(k)=a_2(k)(k+2)^{1-m}k^{1-n}-(k^2-1)B(k).$$ Since $$\lim_{k\to \infty} \frac{s_{2, k}}{k}<\infty\implies \lim_{k\to \infty} \frac{s'_{2, k}}{k}<\infty,$$ we must have $$\deg A-\deg B\leq 3.$$ Since $$\deg B\leq 2\implies \deg A\leq 5.$$ Furthermore, we have already observed that $$A(-6)=A(-5)=A(-4)=A(0)=A(1)=A(2)=0.$$ This implies $A=0$ hence $s'_{2, k}=0$ as claimed. 

The cases $g=3$ and $g=4$ have to be considered separately. First, when $g=4$ we obtain $$m\leq 1,\,\, n\leq 2$$ and $B(k)^2$ divides the polynomial $\widehat v_4(k+2)$. By direct calculation via the Verlinde formula we find $$\widehat v_4(x)=\frac{2x^6+21x^4+168x^2-191}{945}.$$ This implies $B=1$, and thus $$s'_{2, k}=\frac{A(k)}{k^2(k+2)}$$ with $$\deg A\leq 4.$$ Since $A=0$ for $6$ different values, it follows as before that $A=0$ hence $s'_{2, k}=0$.  

When $g=3$, the Verlinde flatness does not give us useful information. In this case, one possible argument is via relative Thaddeus flips, for which we refer the reader to the preprint \cite{FMP}. Along these lines, although we do not explicitly show the details here, the genus $3$ slope formula was in fact checked by direct calculation.
\qed

\section*{Part II: Representation-theoretic methods}

\section {The slope of the Verlinde bundles via conformal blocks}
\label{cfbl}
We derive here the Main Formula \eqref{mainformula} using results in the extensive literature on conformal blocks. In particular, the central statement of \cite{T} is used in an essential way. The derivation is by direct comparison of the bundle ${\mathbb V}_{r, k}^{+}$ of generalized theta functions with the bundle of covacua $$\mathcal B_{r, k}\to \mathcal M_{g,1}$$ defined using the representation theory of the affine Lie algebra $\widehat {\mathfrak {sl}}_r.$ Over pointed curves $(C, p)$, the fibers of the dual bundle $\mathcal B_{r, k}^{\vee}$ give the spaces of generalized theta functions $$H^0(\s\u_C(r, \mathcal O), \left(\Theta_{r}^+\right)^k).$$ Globally, the identification $\mathcal B_{r,k}^{\vee}  \simeq {\mathbb V}_{r,k}^{+}$ will be shown below to hold only up to a twist. The explicit identification of the twist and formula \eqref{mainformula} will be deduced together. 

\subsection {The bundles of covacua}\label{covacua} For a self-contained presentation, we start by reviewing briefly the definition of $\mathcal B_{r,k}.$ Fix a smooth pointed curve $(C, p)$, and write $K$ for the field of fractions of the completed local ring $\mathcal O=\widehat {\mathcal O}_{C, p}.$ For 
notational simplicity, we set $$\mathfrak g=\mathfrak {sl}_r,$$ and write $(|)$ for the suitably normalized Killing form. The loop algebra is the central extension $$\widehat{L\mathfrak g}=\mathfrak g\otimes K\oplus \mathbb C\cdot c$$ of $\mathfrak g\otimes K$, endowed with the bracket $$[X\otimes f, Y\otimes g]=[X, Y]\otimes fg+(X|Y)\cdot \text{Res} \left(g \,df\right)\cdot c.$$ 
Two natural subalgebras of the loop algebra $\widehat{L{\mathfrak g}}$ play a role: $$\widehat{L^+\mathfrak g}=\mathfrak g \otimes \mathcal O\oplus \mathbb C\cdot c \hookrightarrow \widehat{L\mathfrak g}$$ and $$L_{C}\,\mathfrak g=\mathfrak g\otimes \mathcal O_{C} (C-p) \hookrightarrow \widehat{L\mathfrak g}.$$
 
For each positive integer $k$, we consider the basic representation $H_k$ of $\widehat {L\mathfrak g}$ at level $k,$ defined as follows. The one-dimensional vector space $\mathbb C$ is viewed as a module over the universal enveloping algebra $U(\widehat {L^+\mathfrak g})$ where the center $c$ acts as multiplication by $k$, and $\mathfrak g$ acts trivially. We set $$V_k=U(\widehat {L\mathfrak g})\otimes_{U(\widehat {L^+\mathfrak g})} \mathbb C.$$ There is a unique maximal $\widehat{L{\mathfrak g}}$-invariant submodule $$V'_k\hookrightarrow V_k.$$ The basic representation is the quotient $$H_k=V_k/V'_k.$$ The finite-dimensional space of covacua for $(C, p)$, dual to the space of conformal blocks, is given in turn as a quotient $$B_{r, k}=H_k/L_{C}\mathfrak g \,\,H_k.$$ When the pointed curve varies, the loop algebra as well as its two natural subalgebras relativize over $\m_{g,1}.$ The above constructions then give rise to the finite-rank vector bundle $$\mathcal B_{r, k}\to \mathcal M_{g, 1},$$ endowed with the projectively flat WZW connection. 

\subsection {Atiyah algebras} The key theorem in \cite {T} uses the language of Atiyah algebras to describe the WZW connection on the bundles $\mathcal B_{r, k}$. We review this now, and refer the reader to \cite {Lo} for a different account.

An Atiyah algebra over a smooth base $S$ is a Lie algebra which sits in an extension $$0\to \mathcal O_S\stackrel{i}{\rightarrow}\mathcal A\stackrel{\pi}{\rightarrow}\mathcal T_S\to 0.$$ If $L\to S$ is a line bundle, then the sheaf of first order differential operators acting on $L$ is an Atiyah algebra $$\mathcal A_L=\text{Diff}^1(L),$$ via the symbol exact sequence.  

We also need an analogue of the sheaf of differential operators acting on tensor powers $L^c$ for all rational numbers $c$, even though these line bundles don't actually make sense. To this end, if $\mathcal A$ is an Atiyah algebra and $c\in \mathbb Q$, then $c\mathcal A$ is by definition the Atiyah algebra $$c\mathcal A=(\mathcal O_S\oplus \mathcal A)/(c, 1)\mathcal O_S$$ sitting canonically in an exact sequence $$0\to \mathcal O_S\to c\mathcal A\to \mathcal T_S\to 0.$$ The sum of two Atiyah algebras $\mathcal A$ and $\mathcal B$ is given by $$\mathcal A+\mathcal B=\mathcal A\times_{\mathcal T_S}\mathcal B/{(i_{\mathcal A}(f), -i_{\mathcal B}(f)) \text{ for } f\in \mathcal O_S.}$$ When $c$ is a positive integer, $c\mathcal A$ coincides with the sum $\mathcal A+\ldots+\mathcal A$, but $c\mathcal A$ is more generally defined for all $c\in \mathbb Q$. In particular, $c\mathcal A_{L}$ makes sense for any $c\in \mathbb Q$ and any line bundle $L\to S$. 

An action of an Atiyah algebra $\mathcal A$ on a vector bundle $\mathcal V$ is understood to enjoy the following properties 
\begin {itemize}
\item [(i)] each section $a$ of $\mathcal A$ acts as a first order differential operator on $\mathcal V$ with symbol given by $\pi(a)\otimes \mathbf 1_{\mathcal V}$; 
\item [(ii)] the image of $1\in \mathcal O_S$ i.e. $i(1)$ acts on $\mathcal V$ via the identity.  
\end {itemize}
It is immediate that the action of an Atiyah algebra on $\mathcal V$ is tantamount to a projectively flat connection in $\mathcal V$. 
Furthermore, if two Atiyah algebras $\mathcal A$ and $\mathcal B$ act on vector bundles $\mathcal V$ and $\mathcal W$ respectively, then the sum $\mathcal A+\mathcal B$ acts on $\mathcal V\otimes \mathcal W$ via $$(a, b)\cdot v\otimes w=av\otimes w+v\otimes bw.$$ 

We will make use of the following:

\begin {lemma} \label{conn}Let $c\in \mathbb Q$ be a rational number and $L\to S$ be a line bundle. If the Atiyah algebra $c\mathcal A_{L}$ acts on a vector bundle $\mathcal V$, then the slope
$\mu(\mathcal V)= \text{\em det}\,\mathcal V/{\text{\em rank}}\,\mathcal V$ is determined by  $$\mu(\mathcal V)=c\,L \ .$$
\end {lemma}

\proof Replacing the pair $(\mathcal V, L)$ by a suitable tensor power we reduce to the case $c\in \mathbb Z$ via the observation preceding the Lemma. Then, we induct on $c$, adding one copy of the Atiyah algebra of $L$ at a time. The base case $c=0$ corresponds to a flat connection in $\mathcal V$. Indeed, the Atiyah algebra of $\mathcal O_S$ splits as $\mathcal O_S\oplus \mathcal T_S$ and an action of this algebra of $\mathcal V$ is equivalent to differential operators $\nabla_X$ for $X\in \mathcal T_S$, such that $$[\nabla_X, \nabla_Y]=\nabla_{[X, Y]},$$ hence to a flat connection. \qed

\vspace{+10pt}

Consider the rational number $$c=\frac{k(r^2-1)}{r+k},$$ which is the charge of the Virasoro algebra acting on the basic level $k$ representation $H_k$ of $\widehat {L\mathfrak g}$. The representation $H_k$ entered the construction of the bundles of covacua $\mathcal B_{r, k}$. The main result of \cite {T} is the fact that the Atiyah algebra $$\frac{c}{2}\mathcal A_{L}$$ acts on the bundle of covacua $\mathcal B_{r, k}$ where $\mathcal A_{L}$ is the Atiyah algebra associated to the determinant of the Hodge bundle $$L=\det \mathbb E.$$ By Lemma \ref{conn},
we deduce the slope  $$\mu(\mathcal B_{r, k})=\frac{k(r^2-1)}{2(r+k)}\lambda.$$ In fact, by the proof of Lemma \ref{conn}, the bundle $$\mathcal B_{r, k}^{2(r+k)}\otimes L^{-k(r^2-1)}$$ is flat.  

\subsection {Identifications and the slope calculation} We now explain how the above calculation implies the Main Formula \eqref{mainformula} via the results of Section 5.7 of \cite {L}. 

Crucially, Laszlo proves that the projectivization of $\mathcal B_{r, k}^{\vee}$ coincides with the projectivization of the bundle $\mathbb V_{r, k}^+$ coming from geometry. In fact, Laszlo shows that for a suitable line bundle $\mathcal L_{r}$ over $$\s\u_{g, 1}(r, \mathcal O)\to \mathcal M_{g,1}$$ we have $$\mathcal B_{r, k}^{\vee}=\pi_{\star}(\mathcal L_{r}^{k}),$$ where fiberwise, over a fixed pointed curve, $\mathcal L_{r}^{k}$ coincides with the usual theta bundle $\left(\Theta_{r}^+\right)^k.$ Hence, $$\mathcal L_r^{k}=\left(\Theta_{r}^+\right)^{k}\otimes \mathcal T_{r, k}$$ for some line for some line bundle twist $\mathcal T_{r, k}\to \mathcal M_{g, 1}$ over the moduli stack. 
At the heart of this identification is the double quotient construction of the moduli space of bundles over a curve $$\s\u_C(r, \mathcal O)=L_C{G}\backslash \,\widehat {LG}\,/\widehat {L^+G}$$ with the theta bundle $\Theta_{r}^+$ being obtained by descent of a natural line bundle $\mathcal Q_r$ from the affine Grassmannian $$\mathcal Q_r\to \widehat {LG}\,/\widehat {L^+G}.$$ Here $\widehat{LG}$ and $\widehat{L^+G}$ are the central extensions of the corresponding loop groups. The construction is then carried out relatively over $\mathcal M_{g, 1}$, such that $\mathcal Q_r^k$ descends to the line bundle $$\mathcal L_r^k\to \s\u_{g, 1}(r, \mathcal O).$$ It follows from here that fiberwise $\mathcal L_{r}^{k}$ coincides with the usual theta bundle $\left(\Theta_{r}^+\right)^k$.

Collecting the above facts, we find that $$\mathcal B_{r, k}^{\vee}=\mathbb V_{r, k}^+\otimes \mathcal T_{r, k}.$$ Therefore $$-\mu(\mathcal B_{r, k})=\mu(\mathbb V_{r, k}^+)+c_1(\mathcal T_{r, k}).$$ Using Lemma \ref{slopecomparison} we conclude that $$-\frac{k(r^2-1)}{2(r+k)}\lambda=\mu(\mathbb V_{r, k})-\frac{kr-1}{2}\lambda+krc_1(\mathcal L)+c_1(\mathcal T_{r, k}).$$ Simplifying, this yields $$\mu(\mathbb V_{r, k})=\frac{r(k^2-1)}{2(r+k)}\lambda -c_1(\mathcal T_{r, k}) - kr c_1(\mathcal L).$$ Now, the left hand side is a multiple of $\lambda$, namely $s_{r, k}\lambda$. The right hand side must be a multiple of $\lambda$ as well. With $$s'_{r, k}=s_{r, k}-\frac{r(k^2-1)}{2(r+k)}.$$ we find that $$s'_{r, k}\lambda=-c_1(\mathcal T_{r, k}) - kr c_1(\mathcal L).$$ This implies that $s'_{r, k}$ must be an integer by comparison with the right hand side, because the Picard group of $\mathcal M_g$ is generated over $\mathbb Z$ by $\lambda$ for $g\geq 2$, see \cite {AC2}. The fact that $s'_{r, k}\in \mathbb Z$ is enough to prove $$s'_{r, k}=0,$$ which is what we need. 

Indeed, as explained in Section \ref{pf}, Grothendieck-Riemann-Roch for the pushforwards giving the Verlinde numbers shows that $$\lim_{k\to \infty}\frac{s'_{r, k}}{k}<\infty.$$ Writing $$s'_{r, k}=a_r(k)/b_r(k)$$ with $\deg a_r(k)\leq \deg b_r(k)+1$, we see by direct calculation that $$\lim_{k\to \infty} s'_{r, k+1}-2s'_{r, k}+s'_{r, k-1}=0.$$ Since the expression in the limit is an integer, it must equal zero. By induction, it follows that $$s'_{r, k}=A_r k+B_r$$ for constants $A_r, B_r$ that may depend on the rank and the genus. Since $$s'_{r, 0}=s'_{r, -2r}=0$$ by the initial condition in Lemma \ref{ic} and by Proposition \ref{rs}, it follows that $A_r=B_r=0$ hence $s'_{r, k}=0$. 

As a consequence, we have now also determined the twist $\mathcal T_{r, k}=\mathcal L^{-kr}.$ Therefore, the bundle of conformal blocks is expressed geometrically as $$\mathcal B_{r, k}^{\vee}=\mathbb V_{r, k}^+\otimes \mathcal L^{-kr}.$$ We remark furthermore that the latter bundle descends to $\mathcal M_g$. To see this, one checks that  $$(\Theta_{r}^+)^{k}\otimes \mathcal L^{-kr}$$ restricts trivially over the fibers of $\s\u_{g, 1}(r, \mathcal O)\to \mathcal \s\u_g(r, \mathcal O).$ This is a straightforward verification. 

\section{Extensions over the boundary} \label{extension} The methods of \cite {T} can be used to find the first Chern class of the bundle of conformal blocks over the compactification $\overline {\mathcal M}_g$. The resulting formula is stated in Theorem \ref {form} below. In particular the first Chern class contains nonzero boundary contributions, contrary to a claim of \cite{S}.

In genus $0$, formulas for the Chern classes of the bundle of conformal blocks were given in \cite {F}, and have been recently brought to simpler form in \cite {Mu}. In higher genus, the expressions we obtain using \cite{T} specialize to the simpler formulas of \cite {Mu}.

As it is necessary to consider parabolics, we begin with some terminology on partitions. We denote by $\p$ the set of Young diagrams with at most $r$ rows and at most $k$ columns. Enumerating the lengths of the rows, we write a diagram $\mu$ as  $$\mu = (\mu^{1}, \ldots, \mu^{r}), \, \, k \geq \mu^{1} \geq \cdots \geq \mu^{r} \geq 0.$$ The partition $\mu$ is viewed as labeling the irreducible representation of the group $SU( r)$ with highest weight $\mu$, which we denote by $V_{\mu}$. Two partitions which differ by the augmentation of the rows by a common number of boxes yield isomorphic representations. We will identify such partitions in $\mathcal P_{r, k}$, writing $\sim$ for the equivalence relation. There is a natural involution $$\p\ni \mu\mapsto\mu^{\star}\in \p$$ where $\mu^{\star}$ is the diagram whose row lengths are $$k\geq k-\mu^r\geq \ldots \geq k-\mu^1\geq 0.$$
Further, to allow for an arbitrary number of markings, we consider multipartitions $$\underline \mu=(\mu_1, \ldots, \mu_n)$$ whose members belong to $\mathcal P_{r, k}/\sim$. Finally, for a single partition $\mu$, we write $$\mathsf w_{\mu}=\frac{1}{2(r+k)}\left(\sum_{i=1}^{r} \mu_i^2-\frac{1}{r}\left(\sum_{i=1}^r \mu_i\right)^2+ \sum_{i=1}^{r} (r-2i+1)\mu_i\right)$$ for the suitably normalized action of the Casimir element  on the representation $V_{\mu}$. 

In this setup, we let $${\mathcal B}_{g, {\underline \mu}}\to \overline {\mathcal M}_{g, n}$$ be the bundle of covacua, obtained analogously to the construction of Section \ref{covacua} using  representations of highest weight $\underline \mu$, see \cite {T}. To simplify notation, we do not indicate dependence on $r, k$ and $n$ explicitly: these can be read off from the multipartition $\underline \mu$. We set $$v_g(\underline \mu)=\text{rank } \mathcal B_{g, \underline {\mu}}$$ to be the parabolic Verlinde number. 

We determine the  first Chern class $c_1(\mathcal B_{g, \underline \mu})$ over $\overline {\mathcal M}_{g, n}$ in terms of the natural generators:
$$\lambda, \Psi_1, \ldots, \Psi_n$$ and the boundary divisors. To fix notation, we write as usual:
\begin {itemize}
\item $\delta_{\text{irr}}$ for the class of the divisor corresponding to irreducible nodal curves;
\item $\delta_{h, A}$ for the boundary divisor corresponding to reducible nodal curves, with one component having genus $h$ and containing the markings of the set $A$.   
\end {itemize}
Note that each subset $A\subset \{1, 2, \ldots, n\}$ determines a splitting $$\underline \mu_A \cup \underline \mu_{A^c}$$ of the multipartition $\underline \mu$ corresponding to the markings in $A$ and in its complement $A^c$. 
Finally, we define the coefficients $$c_{\text{irr}}=\sum_{\nu\in \mathcal P_{r, k}/\sim} \mathsf w_{\nu} \cdot \frac {v_{g-1}(\underline \mu, \nu, \nu^{\star})}{v_g(\underline \mu)}$$ and $$c_{h, A}=\sum_{\nu\in \mathcal P_{r, k}/\sim} \mathsf w_{\nu} \cdot \frac{v_h(\underline \mu_A, \nu) \cdot v_{g-h}(\underline \mu_{A^c}, \nu^{\star})}{v_g(\underline \mu)}.$$

\begin {theorem} \label{form} Over $\overline {\mathcal M}_{g, n}$ the slope of the bundle of covacua is \begin{equation}\label{slopepar}
\text{\em slope}(\mathcal B_{g, \underline \mu})=\frac{k(r^2-1)}{2(r+k)} \lambda + \sum_{i=1}^n {\mathsf w}_{\mu_i} \Psi_i - c_{\text{irr}} \delta_{\text{irr}} - \sum_{h, A} c_{h, A} \delta_{h, A}.\end{equation}
\end {theorem} 
\noindent In the formula, the repetition $\delta_{h, A}=\delta_{g-h, A^{c}}$ is not allowed, so that each  divisor appears only once. 

\proof The formula written above is correct over the open stratum ${\mathcal M}_{g, n}$. Indeed, the main theorem of \cite {T}, used in the presence of parabolics, shows that the bundle of covacua $$\mathcal B_{g, \underline \mu}\to \mathcal M_{g, n}$$ admits an action of the Atiyah algebra $$\frac{k(r^2-1)}{2(r+k)} \mathcal A_{L}+\sum_{i=1}^{n} \mathsf w_{\mu_i} \mathcal A_{\mathcal L_i}.$$ As before $$L=\det \mathbb E$$ is the determinant of the Hodge bundle and the $\mathcal L_i$ denote the cotangent lines over $\mathcal M_{g, n}$. Therefore, by Lemma \ref{conn}, we have $$\text{slope}(\mathcal B_{g, \underline \mu})=\frac{k(r^2-1)}{2(r+k)} \lambda + \sum_{i=1}^n {\mathsf w}_{\mu_i} \Psi_i$$ over $\mathcal {M}_{g, n}$. 

It remains to confirm that the boundary corrections take the form stated above. Since the derivation is identical for all boundary divisors, let us only find the coefficient of $\delta_{\text{irr}}$. To this end, observe the natural map $$\xi: \mathcal M_{g-1, n+2}\to \overline {\mathcal M}_{g, n}$$ whose image is contained in the divisor $\delta_{\text{irr}}$. The map is obtained by gluing together the last two markings which we denote $\bullet$ and $\star$. We pull back \eqref{slopepar} under $\xi$. For the left hand side, we use the fusion rules of \cite {TUY}: $$\xi^{\star} \mathcal B_{g, \underline \mu}=\bigoplus_{\nu\in \p/\sim} \mathcal B_{g-1, \underline \mu, \nu, \nu^{\star}}.$$ Thus, the left hand side becomes \begin{eqnarray*}\sum_{\nu\in \p/\sim} &{}&\frac {v_{g-1}(\underline \mu, \nu, \nu^{\star})}{v_g(\underline \mu)} \cdot \text{slope} (\mathcal B_{g-1, \underline \mu, \nu, \nu^{\star}})\\&=&\sum_{\nu\in \p/\sim}  \frac {v_{g-1}(\underline \mu, \nu, \nu^{\star})}{v_g(\underline \mu)} \cdot \left(\frac{k(r^2-1)}{2(r+k)} \lambda + \sum_{i=1}^{n} \mathsf w_{\mu_i} \Psi_i+ \mathsf w_{\nu} \Psi_{\bullet}+\mathsf w_{\nu^{\star}} \Psi_{\star}\right)\\&=& \frac{k(r^2-1)}{2(r+k)} \lambda + \sum_{i=1}^{n} \mathsf w_{\mu_i} \Psi_i+ \sum_{\nu\in \p/\sim}  \frac {v_{g-1}(\underline \mu, \nu, \nu^{\star})}{v_g(\underline \mu)} \cdot \left(\mathsf w_{\nu} \Psi_{\bullet}+\mathsf w_{\nu^{\star}} \Psi_{\star}\right).\end{eqnarray*} The fusion rules have been used in the third line to compare the ranks of the Verlinde bundles.
For the right hand side, we record the following well-known formulas \cite {AC}:
\begin {itemize}
\item [(i)] $\xi^{\star}\lambda=\lambda$;
\item [(ii)] $\xi^{\star} \Psi_i=\Psi_i$ for $1\leq i\leq n$;
\item [(iii)] $\xi^{\star} \delta_{\text{irr}}=-\Psi_{\bullet}-\Psi_{\star}$;
\item [(iv)] $\xi^{\star} \delta_{h, A}=0$. 
\end {itemize}
These yield the following expression for the right hand side of \eqref {slopepar}: $$\frac{k(r^2-1)}{2(r+k)}\lambda+  \sum_{i=1}^{n} \mathsf w_{\mu_i} \Psi_i- c_{\text{irr}} (-\Psi_{\bullet}-\Psi_{\star}).$$ For $g-1\geq 2$, $\Psi_{\star}$ and $\Psi_{\bullet}$ are independent in the Picard group of ${\mathcal M}_{g-1, n+2}$, see \cite {AC2}, hence we can identify their coefficient $c_{\text{irr}}$ uniquely to the formula claimed above. The case of the other boundary corrections is entirely similar. \qed

\begin {remark} The low genus case $g\leq 2$ not covered by the above argument can be established
by the following approach.
Once a correct formula for the Chern class has been 
proposed, 
a proof can be obtained by induction on the genus and number of markings. 
Indeed, with some diligent bookkeeping, it can be seen that the expression of the 
Theorem restricts to the boundary divisors compatibly with the fusion rules in \cite {TUY}. 
To finish the argument, we 
invoke the Hodge theoretic result of Arbarello-Cornalba \cite {AC} stating
 the boundary restriction map $$H^2(\overline {\mathcal M}_{g, n})\to H^2(\overline {\mathcal M}_{g-1, n+2})\bigoplus_{h, A}H^2(\overline{\mathcal M}_{h, A\cup\{\bullet\}}\times \overline{\mathcal M}_{g-h, A^c\cup\{\star\}})$$ is injective, with the exception of the particular values $(g, n)=(0, 4), (0, 5), (1, 1), (1, 2)$, which may be checked by hand. 

In fact, the slope expression of the Theorem is certainly correct in the first three cases by \cite {F}, \cite {Mu}. When $(g, n)=(1, 2)$, we already know from \cite {T} that the slope takes the form $$\text{slope}(\mathcal B_{\mu_1, \mu_2})=\frac{k(r^2-1)}{2(r+k)} \lambda + {\mathsf w}_{\mu_1} \Psi_1+{\mathsf w}_{\mu_2}\Psi_2 - c_{\text{irr}} \delta_{\text{irr}}-c \Delta,$$ where $\delta_{\text{irr}}$ and $\Delta$ are the two boundary divisors in $\overline {\mathcal M}_{1, 2}$. The coefficients $c_{\text{irr}}$ and $c$ are determined uniquely in the form stated in the Theorem by restricting $\mathcal B_{\mu_1, \mu_2}$ to the two boundary divisors $\delta_{\text{irr}}$ and $\Delta$ (and not only to their interiors as was done above) via the fusion rules. The verification is not difficult for the particular case $(1,2)$. 

\end {remark}

\end {document}